\documentclass[11pt,reqno]{amsart}
  \usepackage{hyperref}
  \usepackage{amsmath,amssymb}
  \usepackage{xcolor}
  \usepackage{lineno}
  \usepackage{amsthm}
  \usepackage{centernot}
  \usepackage{mathabx}
  \theoremstyle{plain}

  \newtheorem{thm}{Theorem}[section]
   \newtheorem{obs}[thm]{Observation}
  
  \newtheorem{cor}[thm]{Corollary}

  \theoremstyle{definition}
  \newtheorem{defn}[thm]{Definition}

  \newtheorem{rem}[thm]{Remark}
    \newtheorem{note}[thm]{Note}
  \numberwithin{equation}{section}

  \newcommand{\interior}[1]{{\kern0pt#1}^{\mathrm{o}}}

  \makeatletter
  \@namedef{subjclassname@2020}{%
  	\textup{2020} Mathematics Subject Classification}
  \makeatother

\begin{document}
  	
  	\title[Generalizations of Chainability and the Hypertopologies]{Generalizations of Chainability and Compactness, and the Hypertopologies}
  	\author[A. Gupta]{Ajit Kumar Gupta$^{{\MakeLowercase{a}} *}$}
  	\address{$^a$Department of Mathematics\\ National Institute of Technology Meghalaya\\ Shillong 793003\\ India}
  	\email{ajitkumar.gupta@nitm.ac.in}
  	\author[S. Mukherjee]{Saikat Mukherjee{$^{\MakeLowercase{b}}$}}
  	\address{$^b$Department of Mathematics\\ National Institute of Technology Meghalaya\\ Shillong 793003\\ India}
  	\email{saikat.mukherjee@nitm.ac.in}
  	$\thanks{*Corresponding author}$
  	
  	
  	\begin{abstract}
  		We study two properties for subsets of a metric  space. One of them is generalization of chainability, finite chainability, and Menger convexity for metric spaces; while the other is a generalization of compactness. We explore the basic results related to these two properties. Further, in the perspective of these properties, we explore relations among the Hausdorff, Vietoris, and locally finite hypertopologies.
  	\end{abstract}
  	
  	\maketitle

\section{Introduction}\label{intro}

Given the metric space $X$, we denote the collection of nonempty closed subsets of $X$ by $C(X)$. Assuming the space $X$ to be totally bounded, compact, complete, Atsuji, separable, connected, locally connected, or continuum, several researchers have studied the hypertopologies and the relations among the hypertopologies on $C(X)$ (see \cite{gb85, gb87, gb92, jc02, cc95, sg17, ag22, jh99, tj08, sm01, sm04, em51}). A great whim for the study of the hypertopologies arose from the people who were/are more involved in optimization, nonlinear analysis, well-posed problems, etc. than in general point set topology. While studying approximation and/or perturbation methods in analysis, one soon needs the topologies on the set of closed sets. This is not only the matter of analysis, but, the hypertopologies have important applications in probability, statistics, mathematical economics,
game theory, etc. Some applications can be seen in \cite{gb91, gbr92, gbr93, rl95}.

 In \cite{gb85}, Gerald Beer studied the Hausdorff metric topology and the relation of it with the Vietoris topology, assuming that the base space $X$ is Atsuji. He found, a metric space is Atsuji if and only if the Hausdorff metric topology is at least as strong as the Vietoris topology. He also found that, the Hausdorff hyperspace is Atsuji if and only if the space $X$ is either compact or uniformly discrete. Since Atsujiness is weaker than both compactness and uniform discreteness; assuming $X$ to be Atsuji, Gupta and Mekherjee, in the continuation of Beer's work, investigated the Hausdorff hyperspaces (see \cite{ag22}). They proved that if $X$ is Atsuji, then there is a certain class of Atsuji subspaces of the Hausdorff hyperspace. Further, in \cite{gb87}, Beer studied the Hausdorff metric topology and the locally finite topology, assuming that the base space is Atsuji. He proved, $X$ is Atsuji if and only if the locally finite topology coincides with the Hausdorff metric topology. He also proved that, given a metric space $X$, the locally finite topology is metrizable if and only if $X$ admits an Atsuji metric. Jain and Kundu studied various correlations among proximal topology, Vietoris topology, Hausdorff metric topology, and locally finite topology, when the completion of $X$ is Atsuji (see \cite{tj08}).
 
 Michael, in his fundamental paper (see \cite{em51}), proved that a metric space is totally bounded if and only if the Vietoris topology is stronger than the Hausdorff metric topology. While Beer et al. (see \cite{gb92}) showed that the following are equivalent: (1) $X$ is totally bounded; (2) the proximal topology is metrizable; (3) the proximal topology coincides with the Hausdorff metric topology; (4) the Wijsman topology coincides with the proximal topology.
 
 Garc\'ia-Ferreira et al. dealt with the Vietoris hyperspace of all nontrivial convergent sequences $S_c(X)\subset C(X)$ of a space $X$. They showed that (see \cite{sg17}) the space $S_c(X)$ is connected if and only if $X$ is connected. They also proved, the local connectedness of $S_c(X)$ is equivalent to the local connectedness of $X$; and the path-wise connectedness of $S_c(X)$ implies the path-wise connectedness of $X$.
 
 Mac\'ias's work, in \cite{sm01, sm04}, is totally devoted on the continuum spaces and their various hyperspaces. He explored the hyperspaces geometrically in various aspects. We recommend the readers to look at these nice articles.

Chainability for metric spaces is weaker than connectedness, while finite chainability is weaker than total boundedness. The Menger convexity is weaker than the convexity in normed spaces. In this manuscript, we consider two properties for subsets of a metric  space. One of these, named as property $Q$, is weaker than chainability, finite chainability, and Menger convexity for metric spaces. And, the other, which is named as property $I$, is weaker than compactness. We define the two notions and derive the basic results related in Section \ref{Prop. Q n I}. In Section \ref{Viet. n Haus.}, we study the Hausdorff hyperspace when the base space has the property $Q$. We also study the relations between upper/lower Vietoris topology and upper/lower Hausdorff metric topology, with respect to the notions $Q$ and $I$. While, in Section \ref{Loc. Fin. n Haus.} we study the relations between upper/lower locally finite topology and upper/lower Hausdorff metric topology. 

\section{Preliminaries}

  Given a point $x$ and a subset $A$ of the metric space $(X,d)$, the number $\inf\{d(x,y):y\in A\}$ is denoted by $d(x,A)$; and for a subset $B\subset X$, the number $\inf\{d(x,y):x\in A, y\in B\}$ is denoted by $d(A,B)$. We write $B_d(a,\epsilon)$ for an open $\epsilon$-ball with center at $a$, and  $B_d(A,\epsilon)$ for the union of open balls $B_d(a,\epsilon)$ whose centers $a$ run over the set $A$; we shall drop the subscript $d$ while dealing with a single metric on a metric space. The closure of $A$ and the set of all limit points of $A$ are denoted by $\overline{A}$ and $A'$, respectively. We write $\tau_1\lor\tau_2$ for the supremum of two topologies $\tau_1$ and $\tau_2$ on a set $X$. The complement of a set $A$ in $X$ is denoted by $A^c$ or $X\setminus A$; and the collection of nonempty subsets of $X$ is denoted by $P(X)$.

A metric space $(X,d)$ is called \textit{Atsuji} if each continuous function $f:X\to \mathbb R$ is uniformly continuous. Atsujiness is stronger than the completeness, for a metric space. A nonempty subset $A$ in a metric space $(X,d)$ is called \textit{Atsuji} if the subspace $(A,d)$ is an Atsuji space. 
In \cite{tj08}, Jain and Kundu state the following: A completion of a metric space $X$ is Atsuji if and only if every sequence $\{x_n\}\subset X$ with $I(x_n)\to 0$ has a Cauchy subsequence, where $I(x)=d(x,X\setminus \{x\})$.

The \textit{Hausdorff metric} $H$ on $C(X)$ is defined as,
$H(A,B)=\max\{\sup\limits_{x\in A}d(x,B)$, $\sup\limits_{x\in B}d(x,A)\}$, $A,B \in C(X)$. The topology induced by $H$ is denoted by $\tau_H$.
Given the subset $\mathcal B\subset C(X)$, we write $H(A,\mathcal B)$ for the number $\inf\{H(A,B):B\in \mathcal B\}$. The open ball $B_H(A,\epsilon)$ in $C(X)$ is the set $\{B\in C(X):H(A,B)<\epsilon\}$.

\section{Notions More General than Chainability and Compactness}\label{Prop. Q n I}

In this section we will define the properties that are weaker than chainability, finite chainability, compactness, etc.; and study the basic results related to them.

An \textit{$\epsilon$-chain of length $n$} between two points $x,y$ in a metric space $(X,d)$ is a finite sequence $a_0=x,a_1,a_2,...,a_n=y$ in $X$ such that $d(a_{i-1},a_i)\leq \epsilon$ for all $i=1,2,...,n$. A metric space $X$ is said to be \textit{chainable} if for each $\epsilon>0$ each pair of points $x,y$ can be joined by an $\epsilon$-chain. A metric space $X$ is said to be \textit{finitely chainable} if for each $\epsilon>0$, there exist finitely many points $p_1,p_2,...,p_j \in X$ and $m\in \mathbb N$ such that each point of $X$ can be joined with some $p_i$, $1\leq i \leq j$, by an $\epsilon$-chain of length $m$. Chainability is a generalization of connectedness, Menger convexity, etc.; while finite chainability is a generalization of total boundedness. 

In \cite{ag23}, Gupta and Mukherjee introduced a property for metric spaces, which is also a generalization of connectedness, Menger convexity, etc.
\begin{defn}
	A metric space $(X,d)$ is said to have the \textit{property $P$} (\cite{ag23})  if, for all $x\in X$ and for all $r$ with $0< r< \sup\{d(x,y):y\in X\}$, there is a $z\in X$ such that $d(x,z)=r$.
\end{defn}

Consider two distinct elements $x,y$ in a metric space $X$ with the property $P$. For each $n\in \mathbb N$, there is $p_n\in X$ such that $d(x,p_n)=d(x,y)/n$, which implies $x\in X'$. And thus, $X= X'$. Similarly, for a chainable metric space $X$, we have $X=X'$. Notice that, for the spaces in which $X=X'$ holds, the set $I_{\epsilon}=\{x\in X:I(x)>\epsilon\}$ is empty for all $\epsilon>0$, where $I(x)=d(x,X\setminus\{x\})$.  By Proposition 3.3 in \cite{sk17}, in a finitely chainable metric space, the set $I_{\epsilon}$ is finite for all $\epsilon>0$. Thus, these spaces have a common property i.e. the set $I_{\epsilon}$ is finite for all $\epsilon>0$. We consider a functional $\gamma:P(X)\to \mathbb R$, defined as $\gamma (A)=\inf\{\epsilon>0: \text{ the set }\{a\in A:d(a,A\setminus \{a\})>\epsilon\} \text{ is finite}\}$, or alternatively $\gamma (A)=\sup\{\epsilon>0: \text{ the set }\{a\in A:d(a,A\setminus \{a\})>\epsilon\} \text{ is infinite}\}$ and $\gamma(A)=0$ if there is no $\epsilon$ for which the set $\{a\in A:d(a,A\setminus \{a\})>\epsilon\}$ is infinite. Using the functional $\gamma$, we define a property for subsets of a metric space.

\begin{defn}
	A subset $A$ in a metric space is said to have the \textit{property $Q$} if $\gamma(A)=0$.
\end{defn}

 It is clear from the above discussion that the property $Q$ for a metric space is a generalization of:
 
 \begin{enumerate}
 	\item Property $P$,
 	\item Chainability, and
 	\item Finite chainability.
 \end{enumerate}


  The space of rationals $\mathbb Q$ has the property $Q$ but it does not have the property $P$. The subspace $(0,1)\cup (2,3)\subset \mathbb R$ has the property $Q$ but it is not chainable. A normed space has the property $Q$ but it is not finitely chainable. 
  
  It can be verified, a subset $A$ has the property $Q$ if and only if the sequence $I_A(x_n)$ converges to $0$ for each sequence $\{x_n\}\subset A$ of distinct terms, where $I_A(x_n)=d(x_n,A\setminus\{x_n\})$. A known result for Atsujiness is that a metric space is Atsuji if and only if each sequence $\{x_n\}$ with $I(x_n)\to 0$ has a convergent subsequence. Therefore we have, a subset  is compact if and only if it is Atsuji and satisfies the property $Q$.
 
 Notice that, for each open set $U\subset X$ containing the set $I_\epsilon$, we have $\inf\{d(a,U^c):a\in I_\epsilon\}>0$. This suggests a property for subsets of a metric space.
 \begin{defn}
 	A nonempty subset $A$ in a metric space $(X,d)$ is said to satisfy the \textit{property I} if, for each open subset $U$ with $A\subset U$, the infimum of the set $\{d(a,U^c):a\in A\}$ is positive.
 \end{defn}
 
 There are cases where the subsets of a metric space have the property $I$. 
 
 \begin{enumerate}
 	\item If $\inf\{d(a,U^c):a\in A\}=0$ for a given pair of open subset $U$ and subset $A$ of a metric space $(X,d)$ with $A\subset U$, then the infimum $0$ can't be attained at any point of $A$. And hence,  in a metric space, each nonempty compact subset satisfies the property $I$. A metric space itself trivially has the property $I$.
 	\item In a metric space, Atsuji subsets with the property $Q$ satisfy the property $I$.
 	\item  Each isolated point in $(C(X),H)$ has the property $I$ in $X$. As we will see in the proof of Theorem \ref{X-Q iff C(X)-Q} that, $H(A,C(X)\setminus \{A\})>\epsilon>0$ implies $I(a)\geq \epsilon$ for all $a\in A$.

 \end{enumerate}

\subsection{Results Related to Properties $Q$ and $I$}\label{Results on Q n I}

Here we discuss the basic results related to the properties $Q$ and $I$. Precisely we will study the conditions by which the transformed states of a subset preserve the notions $Q$ and $I$.

An Atsuji subset in a metric space may not have the property $I$, based on this fact we have the following observation.
\begin{obs}
	The image of an Atsuji space under a uniformly continuous map need not be Atsuji; even need not be complete. \\ For, given a nonempty subset $A$ in a metric space $X$, and an open subset $U$ with $A\subset U$, it is verified: $\inf\{d(a,U^c):a\in A\}$ is positive if and only if there is an $\epsilon>0$ such that $B(A,\epsilon)\subset U$. Now consider the space $X=\mathbb R^2$ endowed with the Euclidean metric $d$.  The subset $A=\{(n,0):n\in \mathbb N\}$ of $X$ is Atsuji, and the open set $U=\bigcup_{n=1}^\infty B((n,0),1/n)$ contains $A$. If the image of $A$ under every  uniformly continuous map valued in a metric space is Atsuji or complete, then the set $R=\{d(a,U^c):a\in A\}$ is complete in $\mathbb R$, and hence $R$ contains its infimum. Since the infimum is positive, there is an $\epsilon>0$ such that $B(A,\epsilon)\subset U$, a contradiction.
\end{obs}

However, a bi-uniformly continuous map preserves the properties Atsujiness, $Q$, and $I$.
\begin{thm}\label{UnifHom-PropI}
	Consider a uniform homeomorphism $f:X\to Y$ between two metric spaces $X$ and $Y$. Then we have the following: 
	
	\begin{enumerate}
		\item $A\subset X$ has the property $I$ if and only if $f(A)$ has the property $I$;
		\item $A\subset X$ has the property $Q$ if and only if $f(A)$ has the property $Q$;
		\item The completion of $X$ is Atsuji if and only if the completion of $Y$ is Atsuji.
		
	\end{enumerate}
\end{thm}
\begin{proof}
	(1) Suppose $(X,d)$ and $(Y,\rho)$ are the given metric spaces, and a subset $A$ of $X$ has the property $I$. Let $V\subset Y$ be an open subset with $f(A)\subset V$. Consider any two sequences $\{b_n\}\subset f(A)$ and $\{v_n\}\subset V^c$. Then, there exist sequences $\{a_n\}\subset A$ and $\{u_n\}\subset [f^{-1}(V)]^c$ such that $f(a_n)=b_n$ and $f(u_n)=v_n$ for all $n$. Since $A\subset f^{-1}(V)$ and $f^{-1}(V)$ is open, so $d(a_n,u_n)\nrightarrow 0$ when $n\to \infty$. Then by the uniform continuity of $f^{-1}$, $\rho(b_n,v_n)\nrightarrow 0$. And hence, $f(A)$ has the property $I$. The same way we can prove the converse part.

	(2) Suppose $A$ has the property $Q$. Consider a sequence $\{y_n\}\subset f(A)$ with distinct terms. So, there is a sequence $\{x_n\}$ in the subspace $(A,d)$ with distinct terms and $I_A(x_n)\to 0$ as $n\to \infty$, where $f(x_n)=y_n$ for all $n$. This implies there is a sequence $\{p_n\}\subset A$ such that $d(x_n,p_n) \to 0$. By the uniform continuity of $f$, we have $\rho (y_n,f(p_n))\to 0$. And thus, $f(A)$ has the property $Q$. Similarly, the converse can be proved.
	
	(3) Suppose a completion of $X$ is Atsuji, and $\{y_n\}\subset Y$ is sequence with $I(y_n)\to 0$. Then, by the uniform continuity of $f^{-1}$, there is a sequence $\{x_n\}\subset X$ with $I(x_n)\to 0$, where $f(x_n)=y_n$ for all $n$. By the Atsuji completion of $X$, the sequence $\{x_n\}$  has a Cauchy subsequence, say $\{x_{n_k}\}_{k=1}^\infty$; and so, by the uniform continuity of $f$, the subsequence $\{y_{n_k}\}_{k=1}^\infty$ is Cauchy. Thus, a completion of $Y$ is Atsuji. The converse can be proved the same way.
\end{proof}

From $(3)$ of the above theorem, we also have: Given a uniform homeomorphism $f:X\to Y$, $X$ is Atsuji if and only if $Y$ is Atsuji. 

A homeomorphism does not preserve the properties $Q$, $I$, Atsuji completion, or Atsujiness. For example, consider the set $N=\{1/n:n\in \mathbb N\}$ equipped with the usual metric $d_1$ of $\mathbb R$ and with the discrete metric $d_2$. $d_1$ is equivalent to $d_2$. The set $\{1/2n:n\in \mathbb N\}$ has the property $I$ in $d_2$, but not in $d_1$. The set $N$ itself has the property $Q$ in $d_1$, but not in $d_2$. The space $(N,d_2)$ is Atsuji, while $(N,d_1)$ is not. Now, consider the homeomorphism $\tan x:(-\pi/2, \pi/2)\to \mathbb R$. The completion of $(-\pi/2,\pi/2)$ is Atsuji, but $\mathbb R$ is not.

Two metrics $d$ and $\rho$ on a set $X$ are said to be \textit{uniformly equivalent} if, for all $\epsilon>0$ there exist $\alpha, \beta >0$ such that for all $x\in X$ we have $B_d(x,\alpha)\subset B_\rho(x,\epsilon)$ and $B_\rho(x,\beta)\subset B_d(x,\epsilon)$.
\begin{thm} \label{PropI-Unif.Equiv.}
	Suppose $d$ and $\rho$ are equivalent metrics on a set $X$. Then, the following are equivalent:
	\begin{enumerate}
		\item $d$ and $\rho$ are uniformly equivalent;
		\item $d$ and $\rho$ induce the same subsets with the property $I$;
		\item $d$ and $\rho$ induce the same closed subsets with the property $I$.
	\end{enumerate}
\end{thm}
\begin{proof}
	We prove the equivalence of $(1)$ and $(2)$ first.
	
	$(1)\implies (2)$. It follows from (1) of Theorem \ref{UnifHom-PropI}.
	
$(2) \implies (1)$. Since a singleton set satisfies the property $I$, therefore, $d$ and $\rho$ are uniformly equivalent if and only if for each $\epsilon>0$ $\exists \alpha, \beta >0$ such that for all $A\in I^d(X)\cap I^{\rho}(X)$ we have $B_d(A,\alpha)\subset B_\rho(A,\epsilon)$ and $B_\rho(A,\beta)\subset B_d(A,\epsilon)$, where $I^d(X)$ denotes the collection of all nonempty subsets of metric space $(X,d)$ with the property $I$. Now, consider $A\in I^d(X)$. Since $B_\rho(A,\epsilon)$ is open in $(X,d)$, there is an $r>0$ such that $B_d(A,r)\subset B_\rho(A,\epsilon)$. And similarly, for an $A\in I^{\rho}(X)$, there is an $s>0$ such that $B_\rho(A,s)\subset B_d(A,\epsilon)$. And hence proved.

Since the singleton sets in a metric space are closed, we can prove the equivalence of $(1)$ and $(3)$, as above.
\end{proof}

 \begin{cor}
      A metrizable space $X$ is compact if and only if the collection of nonempty closed subsets with the property $I$ remains the same for all the admissible metrics.
 \end{cor}

\begin{thm}\label{AAbarQ}
	If a nonempty subset $A$ in a metric space $X$ has the property $Q$, then so $\overline A$ has. The converse is true if $X'$ is compact.
\end{thm}

\begin{proof}
	Let $A$ have the property $Q$. For $x\in A'$, we have $d(x,\overline A\setminus \{x\})=0$; and $\{x\in A:d(x,\overline A\setminus \{x\})>\epsilon\}\subset \{x\in A:d(x, A\setminus \{x\})>\epsilon\}$ for all $\epsilon>0$. Hence, $\overline A$ has the property $Q$.
	
	For the converse, it is given $X'$ is compact. Suppose $A$ does not have the property $Q$. Then, for some $\epsilon>0$, there is a sequence $\{x_n\}$ in $A$ satisfying $d(x_m,x_n)>\epsilon~ \forall m\neq n$ and $B(x_n,\epsilon)\cap A= \{x_n\}$ for all $n$. None of the $x_n$'s is in $A'$. We claim that there is an $\epsilon'>0$ such that $B(x_n,\epsilon')\cap A'=\emptyset$ for all $n$. If not so, then for all $p\in \mathbb N$ there are $x_{n_p}\in \{x_n:n\in \mathbb N\}$ and $y_p\in A'$ such that $y_p\in B(x_{n_p},1/p)\cap A'$, which implies $d(x_{n_p},y_p)\to 0$ as $p\to \infty$. Since $A'$ is compact, the sequences $\{y_p\}$ and $\{x_{n_p}\}$ have a common limit point, which is a contradiction to $d(x_m,x_n)>\epsilon~ \forall m\neq n$. Now, consider $\epsilon_1=\min \{\epsilon,\epsilon'\}$. Then, $d(x_n,\overline A\setminus \{x_n\})>\frac{\epsilon_1}{2} ~ \forall n$, which implies $\overline A$ does not have the property $Q$. Hence proved.
\end{proof}

\begin{thm}
	If each nonempty closed subset of $X$ has the property $I$, then $X$ is complete.
\end{thm}
\begin{proof}
	If possible, suppose there is a Cauchy sequence $\{x_n\}$ in $X$ which is not convergent. Therefore, the sets $S=\{x_{2n-1}: n\in \mathbb N \}$ and $T=\{x_{2n}: n\in \mathbb N \}$ are disjoint and closed. Now, since a metric space is normal, there are disjoint open sets $U, V$ with $S\subset U$ and $T\subset V$. By the property $I$ of both $S$ and $T$, there are $\epsilon, \epsilon'>0$ such that $B(S,\epsilon)\subset U$ and $B(T,\epsilon')\subset V$. This implies $\{x_n\}$ is not a Cauchy sequence, a contradiction.
\end{proof}

\section{Hausdorff Metric and Vietoris Topologies}\label{Viet. n Haus.}

There are the notions that both a metric space and its Hausdorff hyperspace share with each other, for example, total boundedness, completeness, compactness, etc. (see \cite{gb93}, p. 87). Here we see that they share the property $Q$ as well.

\begin{thm}\label{X-Q iff C(X)-Q}
	A metric space $X$ has the property $Q$ if and only if the hyperspace $(C(X),H)$ has the property $Q$.
\end{thm}
\begin{proof}
	If possible, suppose $C(X)$ does not have the property $Q$. This implies, for some $\epsilon>0$ the collection $\mathcal S=\{A\in C(X):\hat I(A)>\epsilon\}$ is infinite, where $\hat I(A)=H(A,C(X)\setminus \{A\})$. Then, the set $\bigcup_{A\in \mathcal S} A$ is infinite. From the inequality $\hat I(A)>\epsilon$, we claim that $I(a)\geq \epsilon$ for all $a\in A$. For, suppose for some $a\in A$, $I(a)< \epsilon$. This implies, $\exists ~x\in X\setminus \{a\}$ with $d(a,x)< \epsilon$. If $x\not\in A$, then $0<H(A,A\cup \{x\})<\epsilon$. If $x\in A$, consider the nonempty closed set $C=A\setminus B(x,\epsilon) \cup \{x\}$. Then we have $0<H(A,C)\leq \epsilon$, because $a\not \in C$. So we get $\hat I(A)\leq \epsilon$, a contradiction. Thus we get, for all $A\in \mathcal S$ and for all $a\in A$, $I(a)\geq\epsilon$; which means $X$ does not have the property $Q$.
	
	Conversely, let $C(X)$ have the property $Q$. This implies, for all $\epsilon>0$ the collection $\{\{x\}\in C(X):\hat I(\{x\})>\epsilon\}$ is finite.  We prove $\hat I(\{x\})\geq I(x)$ for all $x\in X$. If $\hat I(\{x\})<I(x)$, then $\exists~ B\neq \{x\}$ in $C(X)$ such that $H(B,\{x\})<I(x)$, which implies $\exists ~b (\neq x)\in B$ with $d(b,x)<I(x)$. So we get $I(x)<I(x)$, a contradiction. Thus, we have $\{\{x\}:I(x)>\epsilon\}\subset \{\{x\}:\hat I(\{x\})>\epsilon\}$ for all $\epsilon>0$, and hence proved.
\end{proof}

By using the technique similar to the proof of above theorem, we acquire the following.

\begin{cor}
	The number of isolated points in a metric space $X$ is finite if and only if the number of isolated points in $(C(X),H)$ is finite.
\end{cor}

\begin{obs}\label{no. iso. pts.} From the above discussions, we observe the following as well:	 
	\begin{enumerate}
		
		\item If the number of isolated points in $X$ is $n$, then the number of isolated points in $(C(X),H)$ is exactly $2^n-1$.
		\item If $A_1, A_2,..., A_n$ are the only isolated points in $(C(X),H)$, then each $A_i$, $i=1,2,...,n$, is a finite set and the set of isolated point of $X$ equals $\bigcup_{i=1}^nA_i$.
		\item If the number of isolated points in $(C(X),H)$ is $n$, then there is some $m\in \mathbb N$ such that $n=2^m-1$. And, the cardinality of the union of these isolated points is $m$.
		\item $X=X'$ if and only if $C(X)=[C(X)]'$.
		
	\end{enumerate}
\end{obs}
Now we explore the relations among the notions $Q$, $I$, and the hypertopologies.

For a nonempty subset $A$ of a metric space $X$, define $A^+=\{C\in C(X):C\subset A\}$ and $A^-=\{C\in C(X):C\cap A \neq \emptyset\}$.

A base for the \textit{upper Vietoris topology} $\tau_{V^+}$ on $C(X)$ is the family of all sets of the form $U^+$, where $U$ is a nonempty open subset of $X$; while a subbase for the \textit{lower Vietoris topology} $\tau_{V^-}$ is the family of all sets of the form $U^-$, where $U$ is a nonempty open subset of $X$. The \textit{Vietoris topology} $\tau_V$ on $C(X)$ is the supremum of $\tau_{V^+}$ and $\tau_{V^-}$. A base $\mathcal B_V$ for $\tau_V$ is the collection of all sets of the form 
$\langle U_1, U_2,..., U_n\rangle=\big\{A\in C(X):A\subset \bigcup_{i=1}^nU_i \text{ and }A\cap U_i \neq \emptyset~ \forall i\big\}$, where $U_i$'s are open subsets of $X$ and $n\in \mathbb N$.

In the \textit{upper Hausdorff metric topology} $\tau_{H^+}$ on $C(X)$, a local base at $A\in C(X)$  is the collection of all sets of the form $B_{H^+}(A,\epsilon)=\{C\in C(X):C\subset B(A,\epsilon)\}$, where $\epsilon >0$; while in the \textit{lower Hausdorff metric topology} $\tau_{H^-}$ on $C(X)$, a local base at $A$ is the collection of all sets of the form $B_{H^-}(A,\epsilon)=\{C\in C(X):A\subset B(C,\epsilon)\}$, where $\epsilon >0$. The \textit{Hausdorff metric topology} $\tau_H$ is the supremum of $\tau_{H^+}$ and $\tau_{H^-}$.

The inclusions $\tau_{H^+}\subset \tau_{V^+}$ and $\tau_{V^-}\subset \tau_{H^-}$ hold for any metric \textbf{}space.

\begin{thm}\label{TvAndThd}
	Let $X$ be a metric space. Then, $\tau_{V^+}\subset \tau_{H^+}$ if and only if each nonempty closed subset of $X$ has the property $I$.
\end{thm}

\begin{proof}
	
	%
	It is sufficient to prove that: Given an open set $O\subset X$,
	$\{U^+: \text{ U is open and }U $ $\subset O\}\subset \tau_{H^+}$ if and only if each $C\in O^+$ has the property $I$. Let $C\in O^+$, and $V\subset X$ be an open set with $C\subset V$. Clearly, $C\in (O\cap V)^+$. Therefore for some $\epsilon>0$, $B_{H^+}(C,\epsilon)\subset (O\cap V)^+$, which implies $B(C,\epsilon)\subset O\cap V$.
	
	For the converse, suppose $U$ is an open set contained in $O$. Let $C\in U^+$. By the property $I$ of $C$, there is an $\epsilon>0$ such that $B_{H^+}(C,\epsilon)\subset U^+$.
\end{proof}

\begin{thm}\label{PropQ-LowerHTop.}
	If the local base $\{B_{H^-}(A,\epsilon):\epsilon>0\}\subset \tau_{V^-}$, then each $C\in A^+$ has the property $Q$.
\end{thm}

\begin{proof}
	If possible, suppose a $C\in C(X)$ with $C\subset A$ does not have the property $Q$. Then, for some $\epsilon>0$, the set $P=\{x\in C:d(x,C\setminus \{x\})>\epsilon\}$ is infinite. Therefore, there is a sequence $\{x_n\}$ in $P$ satisfying $d(x_m,x_n)> \epsilon~\forall m\neq n$. Since $B_{H^-}(A,\epsilon/2)\in \tau_{V^-}$, so there are open sets $V_1, V_2, ..., V_k$ in $X$ such that $A\in V_1^-\cap V_2^-\cap ... \cap V_k^-\subset B_{H^-}(A,\epsilon/2)$. Consider $a_i\in A\cap V_i$. Then, $\{x_n:n\in \mathbb N\}\subset A\subset B(\{a_1,a_2,...,a_k\}, \epsilon/2)$, which implies at least two terms, say $x_i,x_j$, of the sequence $\{x_n\}$ are in $B(a_r,\epsilon/2)$ for some $1\leq r\leq k$. And therefore, we have $d(x_i,x_j)<\epsilon$, a contradiction.
\end{proof}

\begin{cor}
	If $\tau_{H^-}\subset \tau_{V^-}$, then each nonempty closed subset of $X$ has the property $Q$.
\end{cor}

\begin{note}
	 The hypertopologies when extended to $N(X)$, the collection of nonempty subsets of a topological space $X$, fail to have the decent separation properties. For example, $(N(X),\tau_V)$ is a $T_1$ space only if $X$ is a discrete space. Here we observe another reason for why we restrict our study of hypertopologies to $C(X)$. The auxiliary result, used in the proof of Theorem \ref{TvAndThd}, can be extended from $C(X)$ to $N(X)$, and hence, the Vietoris topology is weaker than the Hausdorff metric topology on $N(X)$ if and only if $X$ is a discrete metric space.
\end{note}

\section{Locally Finite and Hausdorff Metric Topologies}\label{Loc. Fin. n Haus.}

The property $I$ induces the relations between upper/lower locally finite topology and upper/lower Hausdorff metric topology, that we will discuss in this section.

For a family $\mathcal A$ of nonempty subsets of $X$, define $\mathcal A^-=\{C\in C(X):C\cap A\neq \emptyset~  \forall A\in \mathcal A\}$.

The \textit{upper locally finite topology} $\tau_{lf^+}$ on $C(X) $ has as subbasis the collecection of all sets of the form $U^+$, where $U$ is an open subset of $X$; while the \textit{lower locally finite topology} $\tau_{lf^-}$ has as subbasis the collection of all sets of the form $\mathcal U^-$, where $\mathcal U$ is a locally finite collection of open subsets of $X$. The \textit{locally finite topology} $\tau_{lf}=\tau_{lf^-}\lor \tau_{lf^+}$.

If $X$ is a metrizable space, then $\tau_{lf}=\sup\{\tau_{H_d}: d \text{ is an admissible metric} $ $\text{on X}\}$, where $H_d$ is the Hausdorff metric induced by the metric $d$, and $\tau_{H_d}$ is the topology induced by $H_d$. Therefore, the inclusion $\tau_H\subset \tau_{lf}$ holds for any metric space.

\begin{thm}\label{lf^-H^-}
	If each nonempty closed subset of a metric space has the property $I$, then $\tau_{lf^-}\subset \tau_{H^-}$.
\end{thm}
\begin{proof}
	We shall prove that $\mathcal U^-\in \tau_{H^-}$ for all locally finite open collection $\mathcal U=\{U_i:i\in \Lambda\}$, where $\Lambda$ is an index set. Consider an $A\in \mathcal U^-$. For each $i$, choose an $a_i\in A\cap U_i$, then, since $\mathcal U$ is locally finite, therefore the set $ B=\{a_i:i\in \Lambda\}$ has no limit point in $X$. Now, for each $i\in \Lambda$, there is an $\epsilon_i>0$ such that the open ball $B(a_i,\epsilon_i)\subset U_i$. We claim that $ \inf\limits_{i\in \Lambda}\epsilon_i=\epsilon >0$. For, if $\inf\limits_{i\in \Lambda}\epsilon_i=0$, then there are sequences $\{a_n\}_{n=1}^\infty \subset B$ and $\{x_n\}_{n=1}^\infty$ with $x_n\in X\setminus U_n$ such that $d(a_n,x_n)\to 0$ as $n\to \infty$, where $\{U_n\}_{n=1}^\infty\subset \mathcal U$ is a sequence. And therefore, there is an open set $O:=\bigcup_{n=1}^\infty B(a_n,d(a_n,x_n)/2)$ which contains the closed set $P=\{a_n:n\in \mathbb N\}$, and $\inf\limits_{n\in \mathbb N}d(a_n,O^c)=0$, a contradiction.
	Now, consider an element $C\in B_{H^-}(A,\epsilon)$, then, for each $i\in \Lambda$ there is $c_i\in C$ such that $d(a_i,c_i)<\epsilon$. This implies, $C\cap U_i\neq \emptyset$ for all $i\in \Lambda$, i.e., $C\in \mathcal U^-$. And hence, $\mathcal U^-\in \tau_{H^-}$.
\end{proof}
\begin{thm}\label{lfH}
	Let $X$ be a metric space. Then, $\tau_{lf}\subset \tau_H$ if and only if each nonempty closed subset of $X$ has the property $I$.
\end{thm}
\begin{proof}
	From $\tau_{lf}\subset \tau_H$, we have $V^+\in \tau_H$ for all $V$ open in $X$. Let $A$ be a nonempty closed subset of $X$ and $U$ be an open subset with $A\subset U$. Since $A\in U^+$, therefore $B_H(A,\epsilon)\subset U^+$ for some $\epsilon>0$. Now, consider a nonempty closed subset $B$ with $B\subset B({A,\epsilon/2})$. This implies, $H(A,A\cup B)\leq \epsilon/2$, and therefore $A\cup B\in B_H(A,\epsilon)$. Since $B_H(A,\epsilon)\subset U^+$, so $B\in U^+$. Thus we have $B(A,\epsilon/2)\subset U$, which means $A$ has the property $I$.
	
	The converse follows from Theorems \ref{lf^-H^-} and \ref{TvAndThd}.
\end{proof}

\begin{rem}
	For an Atsuji space $X$, the set $X'$ is necessarily compact. Therefore, using Theorem $2.2$ in \cite{gb87} and Theorem \ref{AAbarQ}, we conclude the following. Suppose each nonempty closed subset in a metric space $X$ has the property $I$.  Then, $A\subset X$ has the property $Q$ if and only if $\overline{A}$ has the property $Q$.
\end{rem}

\begin{cor}
	If each nonempty closed subset in a metric space has the property $I$, then, $\tau_H=\tau_{lf}=\tau_{H^+} \lor \tau_{lf^-}$.  
\end{cor}

\begin{proof}
	The proof follows from Theorems \ref{lfH} and \ref{TvAndThd}.
\end{proof}

  \end{document}